\documentclass[11pt,reqno]{amsart}
\usepackage{fullpage,enumerate}
\usepackage{amsfonts}
\usepackage[all]{xy}
\usepackage{amssymb}
\usepackage{comment}
\usepackage{times}
\usepackage{graphicx}
\usepackage{amsmath,bm}
\usepackage[hidelinks]{hyperref}
\usepackage{tikz-cd}
\usepackage{tikz}
\usepackage{appendix}
\usepackage{color}
\usepackage{float}
\usepackage[margin=1.3in,marginparwidth=1.1in,marginparsep=0.1in]{geometry}
\usepackage{marginnote,color}
\usepackage{import}
\usepackage[normalem]{ulem}
\usepackage{graphicx}

\usepackage{mathrsfs}
\usepackage{geometry}
\usepackage{epsfig,amscd,amsthm,mathtools,fourier}
\numberwithin{equation}{section}

\newtheorem{theorem}{Theorem}[section]
\newtheorem{corollary}[theorem]{Corollary}
\newtheorem{lemma}[theorem]{Lemma}
\newtheorem{proposition}[theorem]{Proposition}
\newtheorem*{Mtheorem}{Main Theorem}

\theoremstyle{definition}

\newtheorem{definition}[theorem]{Definition}

\theoremstyle{remark}
\newtheorem{remark}[theorem]{Remark}

\newcommand{\irr}{\mathrm{irr}}

\newcommand{\RR}{{\mathbb R}}
\newcommand{\ZZ}{{\mathbb Z}}
\newcommand{\PP}{{\mathbb P}}

\newcommand{\TT}{{\mathbb T}}

\newcommand{\CO}{{\mathcal O}}

\newcommand{\CK}{{\mathcal K}}
\newcommand{\CM}{{\mathcal M}}

\newcommand{\CF}{{\mathcal F}}

\newcommand{\CV}{{\mathcal V}}
\newcommand{\CE}{{\mathcal E}}

\newcommand{\CY}{{\mathcal Y}}

\newcommand{\CZ}{{\mathcal Z}}

\newcommand{\CC}{{\mathcal C}}

\newcommand{\oCE}{{\overline{\mathcal E}}}
\newcommand{\oCM}{{\overline{\mathcal M}}}

\newcommand{\Ker}{{\mathrm{Ker}}}
\newcommand{\Coker}{{\mathrm{Coker}}}
\newcommand{\rank}{{\mathrm{rank}}}

\newcommand{\diag}{{\mathrm{diag}}}

\newcommand{\ConH}{{\mathrm{ConvexHull}}}

\newcommand{\bp}{{\mathbf{p}}}
\newcommand{\bx}{{\mathbf{x}}}
\newcommand{\by}{{\mathbf{y}}}
\newcommand{\bz}{{\mathbf{z}}}

\newcommand{\bR}{{\mathbf{R}}}
\newcommand{\bC}{{\mathbf{C}}}

\newcommand{\sO}{{\mathit{0\,}}}

\newcommand{\rvline}{\hspace*{-\arraycolsep}\vline\hspace*{-\arraycolsep}}

\def\:{\colon}

\def\lindn{\left[\mathbb Z^2 : N_0\right]}
\def\bb1{\mathbb 1}
\def\GL{{\rm GL}}
\def\G{{\rm G}}
\def\M{{\rm M}}
\def\lw{{\rm lw}}

\title{A note on elliptic curves on toric surfaces}
\author{Michael M. Barash and Ilya Tyomkin}

\thanks{MB was partially supported by the Israel Science Foundation (grant No. 824/18)}

\address{Department of Mathematics\\
	Ben-Gurion University of the Negev\\P.O.Box 653 \\Be'er Sheva\\ 8410501\\  Israel} \email{barash@post.bgu.ac.il, tyomkin@math.bgu.ac.il}
    
\begin{document}
\tikzcdset{arrow style=tikz, diagrams={>=stealth}}

\begin{abstract}
In this paper, we study the Severi varieties parametrizing integral curves of geometric genus one on polarized toric surfaces in characteristic zero and describe their irreducible components. We show that the irreducible components are in natural bijection with certain affine sublattices of the lattice of characters of the toric surface. The sublattices are described explicitly in terms of the polygon defining the polarization of the toric surface.
\end{abstract}
	
\maketitle
	

\section{Introduction}\label{sec:introduction}

Let $\Delta\subset \RR^2$ be a lattice polygon and $(X_\Delta,\CO(\Delta))$ the corresponding polarized toric surface. Recall that for $g\ge 0$, the {\em Severi variety} $V_{g,\Delta}^\irr$ is the locus of integral curves of geometric genus $g$ in the linear system $|\CO(\Delta)|$ that contain no singular points of $X_\Delta$. It is well known that Severi varieties are quasi-projective and equidimensional of dimension 
\[\dim\left(V_{g,\Delta}^\irr\right)=-K_{X_\Delta}\cdot \CO(\Delta)+g-1=|\partial\Delta\cap M|+g-1,\] 
see, e.g., \cite[Lemma~2.6]{CHT23} and \cite[Theorem~1]{KS13}. We refer the reader to \cite{Tyo24} for a survey on the current state of the art in the geometry of Severi varieties on toric surfaces that will appear as an appendix to the second edition of \cite{GLS07}. An interested reader can also find ample literature on Severi varieties on non-toric surfaces, see, e.g., \cite{CC99, Tes09, CFGK17, Z22, BL23, CGY23, CDGK23}. 

The main result of the current paper is the following theorem calculating the number of irreducible components of the genus-one Severi varieties on polarized toric surfaces. The way we prove our theorem shows not only the equality of numbers but exhibits a natural correspondence between the irreducible components of the Severi variety and the affine sublattices of $\ZZ^2$ satisfying the conditions of the theorem.

\begin{Mtheorem}
    Let $\Delta\subset \RR^2$ be a convex lattice polygon and $V_{1,\Delta}^\irr$ the Severi variety of genus-one curves in the linear system $|\CO(\Delta)|$ on the toric surface $X_\Delta$ over an algebraically closed field of characteristic zero. Then the number of irreducible components of $V_{1,\Delta}^\irr$ is equal to the number of affine sublattices $M\subseteq\ZZ^2$ for which the following two conditions hold: 
    \begin{enumerate}
        \item $\partial\Delta\cap M=\partial\Delta\cap\ZZ^2$, and
        \item $|\Delta^\circ\cap M|\geq 1$,
    \end{enumerate}
    where $\partial\Delta$ and $\Delta^\circ$ denote the boundary and the interior of $\Delta$, respectively.
\end{Mtheorem}

In genus zero, the Severi varieties are known to be irreducible; see, e.g., \cite[Proposition~4.1]{Tyo07}. For higher genera, the irreducibility of the Severi varieties depends both on the polarized surface and the genus. For example, in the case of the complex projective plane, all Severi varieties are irreducible by the celebrated theorem of Harris \cite{Har86}, and similarly for complex Hirzebruch surfaces by \cite{Tyo07}. More generally, the irreducibility of the Severi varieties for a rich class of polarized toric surfaces (including all classical surfaces such as projective plane, Hirzebruch surfaces, and toric del Pezzo surfaces) has been recently established in arbitrary characteristic by Christ, He, and the second author in \cite{CHT25}. 

On the other hand, there are examples of polarized toric surfaces admitting reducible Severi varieties; see \cite{LT23, Tyo14}. In particular, Lang and the second author established a lower bound on the number of irreducible components of the Severi varieties. Namely, \cite[Theorem~A]{LT23} asserts that for $g\ge 1$, the number of irreducible components of the genus-$g$ Severi variety is bounded from below by the number of affine sublattices $M\subseteq\ZZ^2$ for which (1) $\partial\Delta\cap M=\partial\Delta\cap\ZZ^2$ and (2) $|\Delta^\circ\cap M|\geq g$. Although, for a general polygon $\Delta$, this is the best known lower bound, there exists a special family of polygons, called {\em kites}, admitting a sharper lower bound, see \cite[Theorem~B]{LT23}. However, for a kite $\Delta$ and genus $g=1$, the two bounds coincide and are equal to the actual number of irreducible components of $V_{1,\Delta}^\irr$ by \cite[Proposition~3.1]{LT23}. The Main Theorem generalizes this result to arbitrary polarized toric surfaces.

Several approaches have been developed to tackle the irreducibility problem of Severi varieties on toric surfaces: The classical approach of Harris \cite{Har86} based on degeneration techniques and monodromy arguments, the approach of Bourqui \cite{Bou16} based on reduction to positive characteristic and counting curves over finite fields, the recent approach of Christ, He, and the second author based on tropical geometry \cite{CHT23, CHT25}, etc. There have been attempts to extend the classical approach of Harris to general toric surfaces, achieving partial results. In particular, Lang \cite{L20} developed tools to compute the monodromy group of rational curves, and Christ, He, and the second author \cite{CHT22} proved degeneration results for polarized toric surfaces associated to the so-called $h$-transverse polygons.

The approach of the current paper is more elementary and is based on the construction of a covering of an open dense subset of the Severi variety, whose irreducible components can be described explicitly and such that the map to the Severi variety induces a bijection between its irreducible components and those of the Severi variety. The desired covering is a torus fibration over a locus in the moduli space of genus-one curves with marked points $\CM_{1,l}$, given by certain equations on $(E;p_1,\dotsc, p_l)\in \CM_{1,l}$ of the form $\CO_E(\sum_iz_ip_i)\simeq \CO_E$. 

\subsection*{Acknowledgments} 

We are grateful to the anonymous reviewers for their careful reading and valuable suggestions. 

\section{Preliminaries}

\subsection*{Conventions and notation}
Throughout the paper, we work over an algebraically closed field $F$ of characteristic zero. For a matrix $A$, we denote the rows and the columns of $A$ by $\bR_i(A)$ and $\bC_j(A)$, respectively. We denote the standard basis of $\RR^l$ by $e_1,\dotsc, e_l$, and set $\bb1:=\sum_ie_i$ to be the vector all of whose entries are $1$. We often use bold font to denote tuples of numbers or points, e.g., $\bx=(x_1,\dotsc, x_l)$ or $\bp=(p_1,\dotsc, p_l)$. 

\subsection{Lattice polygons}

Let $M\simeq\ZZ^2$ be a lattice of rank two, $N$ the dual lattice, $M_\RR:=M\otimes_\ZZ\RR\simeq\RR^2$ the vector space generated by $M$, and $\Delta\subset M_\RR$ a lattice polygon. Recall that the {\em lattice width} of $\Delta$ is defined to be the minimal length of the projections $n(\Delta)$ for $n\in N$, and is denoted by $\lw_M(\Delta)$. Notice that if $M'\subseteq M$ is a sublattice containing the vertices of $\Delta$, then $M_\RR=M'_\RR$ and $\Delta$ is a lattice polygon with respect to $M'$ too, but its lattice width may be different. Since in what follows we will consider several lattices with respect to which a given polygon is integral, we include the lattice $M$ in the notation above.

\begin{definition}
    Two lattice polygons $\Delta\subset M_\RR$ and $\Delta'\subset M'_\RR$ are called {\em equivalent}, $(\Delta,M)\sim (\Delta',M')$, if there exists an affine integral isomorphism $\xi\: M\to M'$ such that $\xi_\RR(\Delta)=\Delta'$.
\end{definition}

Plainly, equivalent lattice polygons have the same lattice width. Next, we recall the well-known {\em Pick's formula}, \cite{Pic99}, that will also be useful in the proofs below. It asserts the following:
    \[2{\rm Area}(\Delta)=2|\Delta^\circ\cap M|+|\partial\Delta\cap M|-2,\]
where the area is computed with respect to the volume form induced by the lattice $M$, and as usual, $\Delta^\circ$ and $\partial\Delta$ denote the interior and the boundary of $\Delta$, respectively. 

Let $\Delta\subset M_\RR$ be a lattice triangle. Recall that $\Delta$ is called {\em primitive} if $|\Delta\cap M|=3$. Denote by $\Delta_d\subset\RR^2$ the standard triangle with vertices $(0,0),(d,0),(0,d)$. It follows from Pick's formula that any primitive triangle has minimal possible area ${\rm Area}(\Delta)=\frac{1}{2}$, i.e., $(\Delta, M)\sim (\Delta_1,\ZZ^2)$. Similarly, if $\Delta$ is a lattice triangle, all of whose sides have integral length $d$, and such that the affine lattice generated by $\partial\Delta\cap M$ is the whole of $M$, then $(\Delta, M)\sim (\Delta_d,\ZZ^2)$, i.e., $\Delta$ is $d$ times a primitive triangle. Notice that in higher dimensions, the analogous statements are no longer true.

\subsection{Toric varieties}
We assume that the reader is familiar with the basics of toric geometry. For the reader's convenience, we remind a few facts and constructions we use in the paper. They can be found in many introductory texts on toric geometry, e.g., \cite{CLS11}. 

Recall that to a convex lattice polygon $\Delta\subset \RR^2$, one can associate a polarized toric surface $(X_\Delta,\CO(\Delta))$ as follows. Consider the cone $\sigma\subset \RR^3$ over the polygon $\Delta\times\{1\}$. It defines a graded monoid $\sigma\cap \ZZ^3$, where the grading is given by the last coordinate. Then the ${\rm Proj}$ of the monoid algebra $S_\sigma:=F[\sigma\cap \ZZ^3]$ is the desired polarized toric surface 
\[\left(X_\Delta,\CO(\Delta)\right):=\left({\rm Proj}(S_\sigma), \CO_{{\rm Proj}(S_\sigma)}(1)\right).\] 
It is well known that $X_\Delta$ is normal and $\CO(\Delta)$ is very ample. Note that the very ampleness of $\CO(\Delta)$ is true only for polygons and is closely related to the fact that primitive lattice triangles have minimal possible area. In higher dimensions, $\CO(\Delta)$ is always ample but may fail to be very ample.

There is a dimension-preserving correspondence between the faces $\emptyset\ne \CF\subseteq \Delta$ and the torus orbits of $X_\Delta$. In particular, $\Delta$ itself corresponds to the two-dimensional orbit $\TT\subset X_\Delta$, the facets (or sides) $\partial_j\Delta\subset\Delta$ to the one-dimensional orbits $O_j$, and the vertices -- to the zero-dimensional orbits. Furthermore, an orbit $O$ belongs to the closure of an orbit $O'$ if and only if there is an inclusion of the corresponding faces $\CF\subseteq \CF'$. 

Any point of the toric surface $X_\Delta$ admits an equivariant affine neighborhood. In particular, the minimal equivariant affine neighborhood of the generic point of $O_j$ is given by $X_j:=\TT\cup O_j$. Its algebra of functions is given by $F\left[m\,|\, m\cdot n_j\ge 0\right]$, where $n_j$ denotes the primitive inner normal to the facet $\partial_j\Delta$ and dot denotes the usual dot product of vectors in $\RR^2$.

Denote by $k$ the number of facets of $\Delta$. For any $1\le j\le k$, let $O_j$ be the one-dimensional orbit corresponding to a facet $\partial_j\Delta\subseteq \Delta$, and set $D_j:=\overline{O}_j$. Then $K_{X_\Delta}=-\sum_{j=1}^kD_j$, and the class group ${\rm Cl}(X_\Delta)$ is generated by the equivariant divisors $D_j$, for $1\le j\le k$. Furthermore, for any $m\in\ZZ^2$, the divisor of a monomial function $x^m$ is given by \[{\rm div}(x^m)=\sum_{j=1}^k(m\cdot n_j)D_j.\]

\subsection{Homogeneous Smith normal form}

Recall that a matrix $A\in \M_{s\times l}(\ZZ)$ of rank $r$ is said to be in {\em Smith normal form} (SNF) if $A$ is diagonal, the non-zero entries on the main diagonal are the first $r$ entries $\alpha_1,\dotsc, \alpha_r$, and $\alpha_1|\dots|\alpha_r$ are positive integers. A classical result of Smith, \cite{Smi61}, asserts that each orbit of the natural action of $\GL_s(\ZZ)\times \GL_l(\ZZ)$ on $\M_{s\times l}(\ZZ)$ contains a unique SNF. Furthermore, for any $k$, the product of the first $k$ diagonal elements of the SNF in the orbit of a matrix $X$ equals the greatest common divisor of the $k$-by-$k$ minors of $X$. The unique SNF $A$ in the orbit of $X$ is called the {\em Smith Normal Form} of the matrix $X$. Recall that the non-zero entries $\alpha_1,\dotsc, \alpha_r$ of the main diagonal of $A$ are called the {\em invariant factors} of $X$.

Consider the space $\M^h_{s\times l}(\ZZ)=\{X\in \M_{s\times l}(\ZZ)\,|\, X\bb1=0\}$, where $\bb1$ denotes the vector, all of whose entries are $1$. In the proof of our main result, it will be convenient to have a theory of Smith normal forms for the action of the group $\GL_s(\ZZ)\times \GL_l^h(\ZZ)$ on $\M^h_{s\times l}(\ZZ)$, where $\GL_l^h(\ZZ)$ denotes the subgroup of invertible matrices $P\in \GL_l(\ZZ)$ satisfying $P\bb1=\bb1$.

\begin{definition}
     Let $A\in \M^h_{s\times l}(\ZZ)$ be a matrix. We say that $A$ is in {\em homogeneous Smith normal form} (HSNF) if the matrix $A'$ obtained from $A$ by erasing the first column is in Smith normal form.
\end{definition}

\begin{proposition}\label{prop:HSNF}
    Each orbit of the $\GL_s(\ZZ)\times \GL_l^h(\ZZ)$ action on $\M^h_{s\times l}(\ZZ)$ contains a unique HSNF. Moreover, the orbits of this action coincide with the orbits of the action of the subgroup $\GL_s(\ZZ)\times \G$, where $\G\le \GL_l^h(\ZZ)$ is the subgroup of matrices for which $\{0\}\times \ZZ^{l-1}\subset  \ZZ^l$ is invariant.
\end{proposition}
\begin{proof}
    Let $X\in \M^h_{s\times l}(\ZZ)$ be a matrix, $X'$ the matrix obtained from $X$ by erasing the first column, and $Q\in \GL_s(\ZZ)$, $P'\in \GL_{l-1}(\ZZ)$ be such that $QX'(P')^{-1}$ is the SNF of $X'$. Set $u:=(I-P')\bb1$ and 
    \[
        P:=\begin{pmatrix}
            1 & 0\\
            u & P'
        \end{pmatrix}
    \]
    Then $P\in\G$, and $QXP^{-1}$ is in homogeneous Smith normal form. It remains to prove the uniqueness of the HSNF in the $\GL_s(\ZZ)\times \GL_l^h(\ZZ)$-orbit of $X$. Let $A$ be a matrix in HSNF belonging to the orbit of $X$. Then, the invariant factors of $X$ equal the invariant factors of $A$. Since the sum of columns of $A$ is zero, its invariant factors are nothing but the non-zero elements of the diagonal matrix $A'$. Thus, $A'$ is uniquely determined by $X$, and hence so is $A$. 
\end{proof}

\begin{definition}
    The HSNF in the orbit of a matrix $X\in \M^h_{s\times l}(\ZZ)$ is called the {\em homogeneous Smith normal form of $X$}.
\end{definition}

\begin{corollary}
    Let $X\in \M^h_{s\times l}(\ZZ)$ be a matrix, and $\alpha_1,\dotsc,\alpha_r$ its invariant factors. Then, the diagonal above the main diagonal of the homogeneous Smith normal form of $X$ is $(\alpha_1,\dotsc, \alpha_r,0,\dotsc, 0)$. 
\end{corollary}

\section{The proof of the main result}\label{sec:proofMT}

As explained in the introduction, the idea is to find a covering of an open dense subset of the Severi variety, whose irreducible components can be described explicitly and such that the map to the Severi variety induces a bijection between its irreducible components and those of the Severi variety. We construct such a covering in several steps by considering specific loci in the moduli spaces of elliptic curves with marked points and torus fibrations over them. 

Let us first fix the notation. Throughout this section, we fix a convex lattice polygon $\Delta\subset\RR^2$ and denote by $k$ the number of facets of $\Delta$. For $1\le j\le k$, we denote by $\partial_j\Delta$, $l_j$, and $n_j$ the facets of $\Delta$, their integral lengths, and primitive inner normals, respectively. We set $l:=\left|\partial\Delta\cap\ZZ^2\right|=\sum_{j=1}^kl_j$, and denote by $A_\Delta\in \M^h_{2\times l}(\ZZ)$ the matrix, whose first $l_1$ columns are $n_1$, next $l_2$ columns are $n_2$, and so on. We denote by $M_0\subseteq \ZZ^2$ the affine sublattice generated by $\partial\Delta\cap\ZZ^2$, and by $N_0\subset\ZZ^2$ the sublattice generated by the columns of $A_\Delta$, i.e., $N_0=\langle n_1,\dotsc, n_k \rangle$. Plainly, $N_0$ is obtained from the linear sublattice associated to $M_0$ by a rotation by $\pi/2$, which induces a bijection between the intermediate lattices $M_0\subseteq M\subseteq \ZZ^2$ and $N_0\subseteq N\subseteq \ZZ^2$.

\subsection{Useful loci in the moduli spaces}\label{subsec:loci}
Let $\CE\to\CM_{1,1}$ and $\sO\:\CM_{1,1}\to\CE$ be the universal elliptic curve and its zero section. Denote by $\CE^l\to\CM_{1,1}$ its $l$-th power over $\CM_{1,1}$.

\begin{lemma}\label{lem:kerdim}
    Let $A\in \M^h_{s\times l}(\ZZ)$ be a matrix of rank $r$. Consider the homomorphism 
    \[
        \begin{tikzcd}
            \phi_A\:\{0\}\times \ZZ^{l-1}\ar[r, hook] & \ZZ^l\arrow[r,"A"] & \ZZ^s
        \end{tikzcd},
    \]
    and the induced homomorphism $\varphi_A\:\{\sO\}\times\CE^{l-1}\to \CE^s$. Let $\CK_A:=\Ker(\varphi_A)$ be its kernel. Then $\CK_A$ is smooth and equidimensional of dimension $l-r$.
\end{lemma}   

\begin{proof}
    By Proposition~\ref{prop:HSNF}, there exist matrices $P\in \GL_l^h(\ZZ)$ and $Q\in \GL_s(\ZZ)$ such that $QAP^{-1}$ is the HSNF of $A$, and $\{0\}\times \ZZ^{l-1}\subset  \ZZ^l$ is invariant under $P$. We have the following commutative diagram
    \[
        \begin{tikzcd}
            \CE^{l-1} \arrow[r, "\cong"]\arrow[d, "\simeq", "P'"'] &\{\sO\}\times \CE^{l-1} \arrow[rr, hook] \arrow[d, "\simeq", "P"']\arrow[rrrr, bend left, "\varphi_A"] & & \CE^l \arrow[rr, "A"] \arrow[d, "\simeq", "P"'] & & \CE^s \arrow[d, "\simeq", "Q"'] \\
            \CE^{l-1} \arrow[r, "\cong"] &\{\sO\}\times \CE^{l-1} \arrow[rr, hook]\arrow[rrrr, bend right, "\varphi_{QAP^{-1}}"] & & \CE^l \arrow[rr, "QAP^{-1}"] & & \CE^s                 
        \end{tikzcd}
    \]
    where $P'$ is the matrix obtained from $P$ by erasing the first column. Let $\alpha_1,\dotsc, \alpha_r$ be the invariant factors of $A$. Then multiplication by $P$ induces an isomorphism 
    \begin{equation}\label{eq:ircompgen}
         \begin{tikzcd}
            \CK_A\arrow[r, "\sim"] & \CK_{QAP^{-1}}\cong\left(\prod_{i=1}^r\CE[\alpha_i]\right)\times_{\CM_{1,1}}\CE^{l-1-r}             
        \end{tikzcd},
    \end{equation}
    where all products are considered over the moduli space $\CM_{1,1}$, and $\CE[\alpha]$ denotes the kernel of the multiplication by $\alpha$. Now, since for any $\alpha$, $\CE[\alpha]\to \CM_{1,1}$ and $\CE\to \CM_{1,1}$ are smooth of relative dimensions zero and one, respectively, it follows that $\CK_A$ is smooth and equidimensional of dimension $l-1-r+\dim(\CM_{1,1})=l-r$ as asserted. 
\end{proof}

\begin{proposition}\label{prop:KDelta}
    The stack $\CK_{\Delta}:=\CK_{A_\Delta}$ is smooth and equidimensional of dimension $l-2$. The irreducible components of $\CK_{\Delta}$ are parametrized by the intermediate lattices $N_0\subseteq N\subseteq \ZZ^2$. Furthermore, no component $\CK_N$ is contained in the union of the diagonals $D\CE^l\subset\CE^l$ unless $\lw_{M_0}(\Delta)=1$, in which case, $\CK_N\subset D\CE^l$ if and only if $N=N_0$.  
\end{proposition}

\begin{proof}
    Since the $n_i$'s generate the sublattice $N_0\subseteq\ZZ^2$ of full rank, it follows that $\rank(A_\Delta)=2$. Thus, $\CK_{\Delta}$ is smooth and equidimensional of dimension $l-2$ by Lemma~\ref{lem:kerdim}. 
    
    To describe the irreducible components, let us compute the invariant factors $\alpha_1,\alpha_2$ of $A_\Delta$. Let $P\in \GL_l^h(\ZZ)$ and $Q\in \GL_2(\ZZ)$ be such that $QA_\Delta P^{-1}$ is the HSNF of $A_\Delta$ and $\{0\}\times \ZZ^{l-1}\subset  \ZZ^l$ is invariant under $P$. Then multiplication by $Q$ induces an isomorphism 
    \[
    \begin{tikzcd}
            \ZZ^2/N_0=\Coker(\phi_{A_\Delta})\arrow[r, "\sim"] & \Coker(\phi_{QA_\Delta P^{-1}})=(\ZZ/\alpha_1\ZZ)\times(\ZZ/\alpha_2\ZZ).
        \end{tikzcd}
    \]
    The group $\ZZ^2/N_0$ is cyclic since $N_0$ contains a primitive integral vector. Therefore, $\alpha_1=1$ and $\alpha_2=\lindn$ by the structure theorem for finite abelian groups. 
    
    The irreducible components of $\CK_\Delta$ can now be described explicitly following the proof of Lemma~\ref{lem:kerdim}. Indeed, in our case, equation~\eqref{eq:ircompgen} reads as 
    \[
    \CK_\Delta\simeq \CK_{QA_\Delta P^{-1}}=\CE[\alpha_2]\times_{\CM_{1,1}}\CE^{l-1-2}=\left(\coprod_{d|\alpha_2}\CY_1[d]\right)\times_{\CM_{1,1}}\CE^{l-3}=\coprod_{d|\lindn}\left(\CY_1[d]\times_{\CM_{1,1}}\CE^{l-3}\right),
    \]
    where $\CY_1[d]$ is the moduli space of elliptic curves with level-$d$ structure. Notice that the products $\CY_1[d]\times_{\CM_{1,1}}\CE^{l-3}$ are irreducible. Indeed, $\CE^{l-3}\to \CM_{1,1}$ is flat with geometrically irreducible fibers, hence so is $\CY_1[d]\times_{\CM_{1,1}}\CE^{l-3}\to \CY_1[d]$. The moduli space $\CY_1[d]$ is irreducible, and the fiber over its generic point is geometrically irreducible. Thus, there exists a single irreducible component $W$ of $\CY_1[d]\times_{\CM_{1,1}}\CE^{l-3}$ containing the generic fiber. By flatness, any irreducible component $W'$ of $\CY_1[d]\times_{\CM_{1,1}}\CE^{l-3}$ dominates $\CY_1[d]$. Thus, the generic point of $W'$ belongs to the generic fiber, which is equal to $W$, and therefore $W=W'$.
    
    For an intermediate lattice $N_0\subseteq N\subseteq \ZZ^2$, set $\CK_N:=P^{-1}\left(\CY_1[d]\times_{\CM_{1,1}}\CE^{l-3}\right)$, where $d:=[N:N_0]$. Then $\CK_N\subseteq\CK_\Delta$ is an irreducible component. Since $\ZZ^2/N_0$ is cyclic, the index $[N:N_0]$ induces a bijection between the intermediate lattices and the positive divisors of $\lindn$. Therefore, the correspondence $N\mapsto \CK_N$ is the desired bijection between the intermediate lattices and the irreducible components. Moreover, $\CK_N\subseteq\CK_\Delta$ is nothing but the locus of tuples $\left([E,\sO;p_i]\right)_{i=1}^l$ such that $p_1=\sO$ and $\sum_{i=1}^lz_ip_i\in (E,\sO)$ is a point of order $\left[N:N_0\right]$, where \[\bz=\frac{1}{\lindn}{\bR}_2(Q)A_\Delta.\] It is straightforward to check that this bijection is independent of the choice of $Q$ and $P$. We leave the details to the reader since we will not use this fact; cf. Remark~\ref{rem:genirrcomp}.
    
    Let us prove the last assertion of the proposition. An irreducible component of $\CK_\Delta$ belongs to the union of the diagonals $D\CE^l$ if and only if it belongs to a large diagonal. Let $1\le i_1<i_2\le l$ be a pair of integers, and $A_{i_1i_2}\in \M^h_{3\times l}(\ZZ)$ the matrix obtained from $A_\Delta$ by adjoining $e_{i_1}-e_{i_2}$ as the third row. Then a component of $\CK_\Delta$ belongs to the $i_1i_2$-diagonal if and only if $\dim(\CK_\Delta)=\dim(\CK_{A_{i_1i_2}})$. By Lemma~\ref{lem:kerdim}, this is equivalent to $\rank(A_{i_1i_2})=2$. Let us show that $\lw_{M_0}(\Delta)=1$ if and only if there exist $1\le i_1<i_2\le l$ such that $\rank(A_{i_1i_2})=2$.
    
    If $\rank(A_{i_1i_2})=2$, then the $l-1$ vectors $\bC_{i_1}(A_\Delta)+\bC_{i_2}(A_\Delta)$ and $\bC_t(A_\Delta)$ for $t\ne i_1,i_2$ are proportional to each other. Let $1\le j_1<j_2\le k$ be such that $\bC_{i_1}(A_\Delta)=n_{j_1}$ and $\bC_{i_2}(A_\Delta)=n_{j_2}$. Then the condition above implies the following: $j_1\ne j_2$, $l_{j_1}=l_{j_2}=1$, and the facets $\partial_j\Delta$ with $j\ne j_1,j_2$ are parallel to each other. In particular, $k\le 4$. Therefore, $\Delta$ has $M_0$-width one with respect to the normal to a facet $\partial_j\Delta$ with $j\ne j_1,j_2$. Vice versa, if $\lw_{M_0}(\Delta)=1$, then $k\le 4$, and there exist $1\le j_1<j_2\le k$ such that  $l_{j_1}=l_{j_2}=1$ and the facets $\partial_j\Delta$ with $j\ne j_1,j_2$ are parallel to each other. Therefore, $\rank(A_{i_1i_2})=2$ where $i_1$ and $i_2$ are such that $\bC_{i_1}(A_\Delta)=n_{j_1}$ and $\bC_{i_2}(A_\Delta)=n_{j_2}$. 
     
    To complete the proof, it remains to show that if $\rank(A_{i_1i_2})=2$, then $\CK_{A_{i_1i_2}}=\CK_{N_0}$. Assume that $\rank(A_{i_1i_2})=2$. We claim that the invariant factors of $A_{i_1i_2}$ are $\beta_1=\beta_2=1$. To see this, pick $i\ne i_1,i_2$. Then the vector $\bC_i(A_\Delta)$ is primitive, and therefore the greatest common divisor of the two-by-two minors in the columns $i$ and $i_1$ of $A_{i_1i_2}$ is $1$. It follows that $\beta_1=\beta_2=1$. Consider now the matrix $\diag(Q,1)A_{i_1i_2}P^{-1}$. It has rank two and contains the matrix $QA_\Delta P^{-1}$ as its top block. Thus, there exist integers $x,y\in\ZZ$ such that  
    \[
        \diag(Q,1)A_{i_1i_2}P^{-1}=
        \begin{pmatrix}
            \begin{matrix}
                -1 & 1 & 0\\
                -\alpha_2 & 0 & \alpha_2\\
                -x-y & x & y
            \end{matrix}
            & \rvline &
            \begin{matrix}
                0 & \dots & 0\\
                0 & \dots & 0\\
                0 & \dots & 0
            \end{matrix}
        \end{pmatrix}
    \]
    Since both invariant factors of $A_{i_1i_2}$ are $1$, it follows that $P(\CK_{A_{i_1i_2}})=\CY_1[1]\times_{\CM_{1,1}}\CE^{l-3}=P(\CK_{N_0})$, which implies $\CK_{A_{i_1i_2}}=\CK_{N_0}$, and we are done.
\end{proof}

\begin{remark}\label{rem:genirrcomp}
    It is possible to describe the irreducible components of $\CK_A$ for an arbitrary matrix $A\in \M^h_{s\times l}(\ZZ)$ of rank $r$. Let $d:=\alpha_1\cdot\dotsc\cdot\alpha_r$ be the product of the invariant factors of $A$, and $\CY_0[d]$ the moduli space of elliptic curves with the full level-$d$ structure. Then the irreducible components of $\CK_A$ are the images of natural morphisms $\CY_0[d]\times_{\CM_{1,1}}\CE^{l-1-r}\to\{\sO\}\times\CE^{l-1}$ associated to the elements of $\left({\rm Tor}^1_\ZZ\left(\ZZ/d\ZZ, \Coker(A)\right)\right)^2$. Furthermore, two such morphisms share the same image if and only if the corresponding elements belong to the same orbit of the right action of ${\rm GL}_2(\ZZ/d\ZZ)$ on $\left({\rm Tor}^1_\ZZ\left(\ZZ/d\ZZ, \Coker(A)\right)\right)^2$.
\end{remark} 

The universal curve $\oCE\to \oCM_{1,1}$ over the compactified moduli space can be identified naturally with the forgetful morphism $\oCM_{1,2}\to \oCM_{1,1}$ that forgets the second marked point and stabilizes the curve. Therefore, the forgetful morphisms $\pi_i\:\oCM_{1,l}\to \oCM_{1,2}$ that forget all marked points but the first and the $(i+1)$-st, induce a natural proper birational morphism from $\oCM_{1,l}$ to the $(l-1)$-th power of $\oCE$ over $\oCM_{1,1}$, which restricts to an open immersion \[\iota\:\CM_{1,l}\to \CE^{l-1}=\{\sO\}\times \CE^{l-1}\subset\CE^l.\] Set $\CZ_\Delta:=\iota^{-1}(\CK_\Delta)\subset\CM_{1,l}$ and $\CZ_N:=\iota^{-1}(\CK_N)\subset\CM_{1,l}$ for the intermediate lattices $N_0\subseteq N\subseteq \ZZ^2$. Plainly, the $\CZ_N$ are either empty or irreducible. The following is an immediate corollary of Proposition~\ref{prop:KDelta}.

\begin{corollary}\label{cor:compCZ}
    In the notation above, the following holds:
    \begin{enumerate}
        \item $\CZ_\Delta$ is smooth and equidimensional of dimension $l-2$;
        \item $\CZ_N=\emptyset$ if and only if\, $N=N_0$ and\, $\lw_{M_0}(\Delta)=1$;
        \item the irreducible components of $\CZ_\Delta$ are the non-empty $\CZ_N$'s.
    \end{enumerate}
\end{corollary}

\begin{remark}\label{rem:expCZ}
    Since for any elliptic curve $(E,\sO)$, there is a natural group isomorphism $E\to {\rm Pic}^0(E)$ given by $p\mapsto \CO_E(p-\sO)$, and since the rows $\bx:=\bR_1(A_\Delta)$ and $\by:=\bR_2(A_\Delta)$ satisfy: $\bx\cdot\bb1=\by\cdot\bb1=0$, the locus $\CZ_\Delta\subset\CM_{1,l}$ can be described more symmetrically as follows
    \[\CZ_\Delta=\left\{(E;\bp)\in\CM_{1,l}\,|\, \CO_E(\bx\bp)\simeq\CO_E(\by\bp)\simeq\CO_E\right\},\]
    where $\bp=(p_1,\dotsc,p_l)$ is the tuple of marked points, $\bx\bp=\sum_{i=1}^l x_ip_i$, and $\by\bp=\sum_{i=1}^l y_ip_i$. Moreover, if $QA_\Delta P^{-1}$ is the HSNF of $A_\Delta$, then by construction, the component $\CZ_N$ is the locus of marked curves $(E;\bp)\in\CZ_\Delta$ for which the order of $\CO_E(\bz\bp)\in {\rm Pic}^0(E)$ is $\left[N:N_0\right]$, where \[\bz=\frac{1}{\lindn}\bR_2(Q)A_\Delta.\] Notice that the vector $\bz$ is invariant under the natural action of the subgroup $\prod_{j=1}^k{\mathfrak S}_{l_j}$ of the full symmetric group ${\mathfrak S}_l$ on $\ZZ^l$ since so are the rows of $A_\Delta$. This observation will be useful in the proof of Corollary~\ref{cor:compcorr}.
\end{remark}

\subsection{The construction of the desired covering}

Let $\pi\:\CC\to \CM_{1,l}$ be the universal curve, and $p_i\:\CM_{1,l}\to \CC$ the sections. Consider the restriction $\CC_{\CZ_\Delta}\to\CZ_\Delta$ of $\CC$ to $\CZ_\Delta$, and the line bundles $\CO_{\CC_{\CZ_\Delta}}(\bx\bp)$ and $\CO_{\CC_{\CZ_\Delta}}(\by\bp)$. By construction of $\CZ_\Delta$, the restrictions of these bundles to the fibers of $\pi$ over $\CZ_\Delta$ are trivial, and therefore $\pi_*\CO_{\CC_{\CZ_\Delta}}(\bx\bp)$ and $\pi_*\CO_{\CC_{\CZ_\Delta}}(\by\bp)$ are invertible sheaves on $\CZ_\Delta$ by Grothendieck's cohomology and base change theorem. 

Consider the total space $\CV$ of the locally free sheaf $\pi_*\left(\CO_{\CC_{\CZ_\Delta}}(\bx\bp)\oplus\CO_{\CC_{\CZ_\Delta}}(\by\bp)\right)$ on $\CZ_\Delta$, and let $\CV_\Delta$ be the complement of the coordinate line bundles in $\CV$. Let $\pi_\Delta\:\CC_\Delta\to \CV_\Delta$ be the pullback of $\pi\:\CC\to \CM_{1,l}$ with respect to  $\CV_\Delta\to \CZ_\Delta\subset\CM_{1,l}$. Since $\CV_\Delta\to\CZ_\Delta$ is flat, surjective, and has irreducible fibers, it follows that the irreducible components of $\CV_\Delta$ are the preimages of those of $\CZ_\Delta$. We denote the preimage of $\CZ_N$ by $\CV_N$. By construction, the points of $\CV_\Delta$ are triples $(q,f_x,f_y)$, where $q\in \CZ_\Delta$ is a point, and $f_x\in H^0(\CC_q,\CO_{\CC_q}(\bx\bp))$ and $f_y\in H^0(\CC_q,\CO_{\CC_q}(\by\bp))$ are a pair of non-zero sections. Therefore, we have a natural rational map $\CC_\Delta\dashrightarrow \TT\subset X_\Delta$ sending $p\in\CC_{(q,f_x,f_y)}$ to $f(p):=\left(f_x(p),f_y(p)\right)$.

\begin{lemma}
    The rational map defined above extends to a morphism $f\:\CC_\Delta\to X_\Delta$, which induces a morphism $\psi\:\CV_\Delta\to |\CO(\Delta)|$ given by $\psi(q,f_x,f_y)=f_*(\CC_q)$. Furthermore, $f^*(\sum_j D_j)=\sum_{i} p_i$, where as usual, $D_j=\overline{O}_j$ is the equivariant divisor corresponding to the facet $\partial_j\Delta$ of $\Delta$.
\end{lemma}
\begin{proof}
    By construction, the map $f$ is defined on the complement of the marked points $\bp$. Furthermore, ${\rm div}\left(f^*(x^m)\right)=\sum_i(m\cdot \bC_i(A_\Delta))p_i$ for any $m\in M$. For an $1\le i\le l$, let $j$ be such that $n_j=\bC_i(A_\Delta)$. Then, along the section $p_i$, the rational map $f$ extends to a morphism to the affine toric subvariety $X_j=\TT\cup O_j$, since the algebra of functions on $X_j$ is given by $F\left[m\,|\, m\cdot n_j\ge 0\right]$. In particular, $f$ extends to a morphism $f\:\CC_\Delta\to X_\Delta$. Furthermore, since for $m\cdot n_j=1$, the restriction of the divisor ${\rm div}(x^m)$ to $X_j=\TT\cup O_j$ is $O_j$, it follows that the multiplicity of $p_i$ in the pullback of the boundary divisor is $1$, i.e., $f^*(\sum_j D_j)=\sum_{i} p_i$. Finally, since $f|_{\CC_q}\:\CC_q\to X_\Delta$ is proper, by the projection formula we have:
    \[f_*(\CC_q)\cdot D_j=\CC_q\cdot f^*(D_j)=l_j=\deg\left(\CO(\Delta)|_{D_j}\right)\] 
    for any $1\le j\le k$; and since the equivariant irreducible divisors $D_j$ generate the class group ${\rm Cl}(X_\Delta)$, it follows that $f_*(\CC_q)\in |\CO(\Delta)|$ for any $(q,f_x,f_y)\in \CV_\Delta$. Therefore, $\psi\:\CV_\Delta\to |\CO(\Delta)|$, as asserted.
\end{proof}

\begin{proposition}\label{prop:proppsi}
    Let $f\:\CC_\Delta\to X_\Delta$ and $\psi\:\CV_\Delta\to |\CO(\Delta)|$ be as in the lemma. Then $V_{\Delta,1}^\irr\subseteq \overline{\psi(\CV_\Delta)}$. Furthermore, the following holds:
    \begin{enumerate}
        \item If the morphism $f|_{\CC_q}$ is birational onto its image, then $\psi(q,f_x,f_y)\in V_{\Delta,1}^\irr$;
        \item For a general point $(q,f_x,f_y)$ in an irreducible component $\CV_N\subseteq \CV_\Delta$, the morphism $f|_{\CC_q}$ is not birational onto its image if and only if $N=N_0$ and $(\Delta,M_0)\sim (\Delta_2,\ZZ^2)$, where $\Delta_2\subset\RR^2$ denotes the triangle with vertices $(0,0),(2,0),(0,2)$.
    \end{enumerate}
\end{proposition}

\begin{proof}
    Pick a general $[C]\in V_{\Delta,1}^\irr$, and let $E:=C^\nu$ be the normalization of $C$ and $f\: E\to X_\Delta$ the natural morphism. By \cite[Proposition~17]{KS13}, $C$ intersects the boundary divisor $\sum_j D_j$ transversally away from the zero-dimensional orbits. Thus, the pullback $f^*(D_j)$ is a reduced divisor of degree $l_j$. Let us label the points in the support of $f^*(\sum_j D_j)$ as $p_i$'s in such a way that $f$ maps the first $l_1$ points $p_1,\dotsc, p_{l_1}$ to $O_1$, the next $l_2$ points to $O_2$, etc. Set $f_x:=f^*(x^{e_1})$ and $f_y:=f^*(x^{e_2})$. Since ${\rm div}(x^m)=\sum_{j=1}^k(m\cdot n_j)D_j$ for any $m\in\ZZ^2$, it follows that \[{\rm div}\left(f^*(x^m)\right)=\sum_{i=1}^l(m\cdot \bC_i(A_\Delta))p_i.\] In particular, ${\rm div}(f_x)=\bx\bp$ and ${\rm div}(f_y)=\by\bp$, and therefore $(E;\bp)\in \CK_\Delta$. Set \[q:=\iota^{-1}(E;\bp)\in\CZ_\Delta.\] Then $(q,f_x,f_y)\in\CV_\Delta$. Thus, $[C]=\psi(q,f_x,f_y)$ and $V_{\Delta,1}^\irr\subseteq \overline{\psi(\CV_\Delta)}$ as asserted.

    (1) Set $[C]:=f_*(\CC_q)\in |\CO(\Delta)|$. If the morphism $f|_{\CC_q}$ is birational onto its image, the curve $C$ is reduced and has geometric genus $1$. Furthermore, it contains no zero-dimensional orbits of $X_\Delta$; in particular, no singular points of $X_\Delta$. Thus, $\psi(q,f_x,f_y)=[C]\in V_{\Delta,1}^\irr$ as asserted.

    (2) If $(\Delta,M_0)\sim (\Delta_2,\ZZ^2)$, then after a change of coordinates in $\ZZ^2$, we may assume without loss of generality that $n_1=e_2$ and $\Delta$ is the triangle with vertices $(0,0),(2,0),(2a,2b)$, where $a$ and $b$ are positive integers satisfying $\gcd\{a,b\}=\gcd\{a-1,b\}=1$. Thus, 
    \[A_\Delta=
        \begin{pmatrix}
            0 & 0 & -b & -b & b & b\\
            1 & 1 & a-1 & a-1 & -a & -a
        \end{pmatrix}.\]
    Set \[Q:=\begin{pmatrix}
            0 & 1\\
            1 & 0
        \end{pmatrix}\]
    and \[P:=\begin{pmatrix}
            1 & 0 & 0 & 0 & 0 & 0\\
            -a & 1 & a-1 & 0 & 1 & 0\\
            0 & 0 & -1 & 1 & 1 & 0\\
            2 & 0 & 0 & -1 & 0 & 0\\
            1 & 0 & 0 & 0 & 1 & -1\\
            0 & 0 & 0 & 0 & 0 & 1
        \end{pmatrix}^{-1}.\]
    Then $P\in \GL_6^h(\ZZ)$ and $QA_\Delta P^{-1}$ is the HSNF of $A_\Delta$ and
    \[\frac{1}{\lindn}\bR_2(Q)A_\Delta= \begin{pmatrix}
                                            0 & 0 & -1 & -1 & 1 & 1
                                        \end{pmatrix}.\]
    Thus, the component $\CZ_{N_0}$ is the locus of marked curves $(E;\bp)\in\CZ_\Delta$ for which $p_3+p_4=p_5+p_6$ by Remark~\ref{rem:expCZ}. Since $\CZ_\Delta\subset\CM_{1,6}$ is given by the equations $p_1+p_2=(1-a)(p_3+p_4)+a(p_5+p_6)$ and $b(p_3+p_4)=b(p_5+p_6)$, it follows that $\CZ_{N_0}\subset\CM_{1,6}$ is given by the following equation:
    \begin{equation}\label{eq:ZN0}
    p_1+p_2=p_3+p_4=p_5+p_6.
    \end{equation} 
    
    By \eqref{eq:ZN0}, for any point $(q,f_x,f_y)\in\CV_{N_0}$, there exists a function $f_0$ on the elliptic curve $\CC_q$ for which ${\rm div}(f_0)=p_5+p_6-p_3-p_4$, and the divisors of $f_x$ and $f_y$ are given by
    \[{\rm div}(f_x)=b{\rm div}(f_0)\quad\text{and}\quad {\rm div}(f_y)=p_1+p_2-p_3-p_4-a{\rm div}(f_0).\]
    Therefore, $f_x$ is a scalar multiple of $f_0^b$, and $f_yf_0^a\in H^0(\CC_q,\CO_{\CC_q}(p_3+p_4))$. The latter space is two-dimensional by the Riemann-Roch Theorem, and it is generated by the functions $1,f_0$. We conclude that $f_x$ and $f_y$ belong to the function space generated by $f_0$, and hence the map $f|_{\CC_q}$ is not birational onto its image.
    
    In the opposite direction, assume that $(q,f_x,f_y)\in \CV_N$ is general and the morphism $f|_{\CC_q}$ is not birational onto its image. Set $d:=\deg(f|_{\CC_q})\ge 2$ and $C:=f(\CC_q)$. Then, for any $1\le j\le k$, the degree $d$ divides $l_j$, and therefore $l_j\ge 2$. Assume first that the curve $C$ is not rational. Then it has geometric genus $1$. The morphism $f|_{\CC_q}$ factors through the normalization $C^\nu$ of $C$, so we obtain a morphism $(\CC_q,p_1)\to (C^\nu,f(p_1))$ between elliptic curves. Since any morphism between abelian varieties is a composition of a translation with a homomorphism, the morphism $(\CC_q,p_1)\to (C^\nu,f(p_1))$ is an isogeny. Since $l_1\ge 2$, there exists $1<i\le l_1$ such that $f(p_i)=f(p_1)$. Therefore, $dp_i=0$ in the group $(\CC_q,p_1)$, and therefore, the irreducible component $\CZ_N$ belongs to the kernel $\CK_A$, where $A$ is the matrix obtained from $A_\Delta$ by adjoining the third row $de_1-de_i$. However, $\rank(A)=3$ since the columns $1,i$, and $l_1+1$ of $A$ are linearly independent. Thus, $\dim(\CZ_N)\le l-3$ by Lemma~\ref{lem:kerdim}, which is a contradiction since $\dim(\CZ_N)=l-2$ by Corollary~\ref{cor:compCZ}. Therefore, $C$ is rational, and its normalization is $\PP^1$.

    Next, let us show that $N=N_0$ and $(\Delta,M_0)\sim (\Delta_d,\ZZ^2)$, where $\Delta_d$ is the triangle with vertices $(0,0),(d,0)$, and $(0,d)$. Since $C^\nu\simeq\PP^1$, the fibers of $\CC_q\to C^\nu$ are linearly equivalent as divisors on $\CC_q$, and the set of marked points $\{p_i\}_{i=1}^l$ splits into a disjoint union of $l/d$ such fibers. In other words, there is a partition $\sqcup_{t=1}^{l/d} I_t$ of $\{1,2,\dotsc, l\}$, such that $\sum_{i\in I_t}p_i=\sum_{i\in I_1}p_i$ in the group $(\CC_q,p_1)$ for any $1\le t\le l/d$. Consider the matrix $A$ with rows $\sum_{i\in I_t}e_i-\sum_{i\in I_1}e_i$ for $t=1,\dotsc, l/d$. By construction, $\CZ_N\subseteq\CK_A$ and $A$ has rank $l/d-1$. Thus, \[l-2=\dim(\CZ_N)\le \dim(\CK_A)=l-l/d+1\]
    by Lemma~\ref{lem:kerdim} and Proposition~\ref{prop:KDelta}, and therefore $l/d\le 3$. Since $d|l_j$ for all $j$, it follows that $k=3$ and $l_1=l_2=l_3=d$, i.e., $(\Delta,M_0)\sim (\Delta_d,\ZZ^2)$ as asserted. To see that $N=N_0$, notice that the two invariant factors of the matrix $A$ are $1$. Let $\hat{A}$ be the matrix obtained from $A_\Delta$ by adding $A$ at the bottom. Then $\CZ_N\subseteq \CK_{\hat{A}}$, and therefore $\rank(\hat{A})=\rank(A_\Delta)=\rank(A)=2$ by Lemma~\ref{lem:kerdim}. Thus, the two invariant factors of $\hat{A}$ are also $1$ since so are the invariant factors of $A$. Therefore $\CZ_N=\CZ_{N_0}$, which implies $N=N_0$.

    Finally, let us show that $d=2$. Assume to the contrary that $d\ge 3$. Let $m=(a,b)\in M$ be such that $m\cdot n_1=0$. Then $f^*(x^m)=f_x^af_y^b$ is a non-constant function defined and not zero at the points $p_1,\dotsc, p_d$. Furthermore, \[{\rm div}\left(f^*(x^m)\right)=(m\cdot n_2)(p_{d+1}+\cdots +p_{2d})+(m\cdot n_2)(p_{2d+1}+\cdots +p_{3d}).\]
    In particular, $f^*(x^m)$ is uniquely determined by the points $p_{d+1},\dotsc, p_{3d}$ up to a multiplicative scalar, and is independent of $p_1,\dotsc,p_d$. Consider the rational map $\theta\:\CC_q\dashrightarrow \CM_{1,l}$, given by \[\theta:\epsilon\mapsto (\CC_q;p_1,p_2+\epsilon,p_3-\epsilon,p_4,\dotsc,p_l),\]
    where the addition and subtraction are, as usual, in the group $(\CC_q,p_1)$. It is defined in a neighborhood of $p_1\in\CC_q$ and $\theta(p_1)=q\in\CM_{1,l}$. Furthermore, the image of $\theta$ belongs to the kernel $\CZ_\Delta=\sqcup_N\CZ_N$, and therefore it belongs to the irreducible component $\CZ_N$. Since $(q,f_x,f_y)\in \CV_N$ is general, it follows that for any $(q',f'_x,f'_y)\in \CV_N$, the morphism $f|_{\CC_{q'}}$ has degree $d$. However, for a general $\epsilon\in\CC_q$, the map $f|_{\CC_{\theta(\epsilon)}}$ separates the points $p_2+\epsilon$ and $p_3-\epsilon$ since so does $f^*(x^m)$, which leads to a contradiction. We conclude that $d=2$.
\end{proof}

\begin{corollary}\label{cor:compcorr}
    The map $\psi$ induces a natural bijection between the set of irreducible components of $V_{\Delta,1}^\irr$ and the set of irreducible components of\; $\CV_\Delta$, unless $(\Delta,M_0)\sim (\Delta_2,\ZZ^2)$, in which case, the bijection is with the subset $\{\CV_N\}_{N\ne N_0}$. 
\end{corollary}
\begin{proof}
    If for a general point $(q,f_x,f_y)\in \CV_N$, the morphism $f|_{\CC_q}$ is not birational onto its image, then $\psi(\CV_N)\subset |\CO(\Delta)|$ is disjoint from the Severi variety. Otherwise, $\psi(q,f_x,f_y)$ belongs to the Severi variety by Proposition~\ref{prop:proppsi} (1), and the map $\psi|_{\CV_N}$ is generically finite. Thus, \[\dim(\psi(\CV_N))=\dim(\CV_N)=\dim(\CZ_N)+2=l=|\partial\Delta\cap\ZZ^2|+1-1=\dim(V_{\Delta,1}^\irr).\]
    It follows from Proposition~\ref{prop:proppsi}, that the map $\psi$ induces a natural surjection $\CV_N\mapsto \overline{\psi(\CV_N)}\cap V_{\Delta,1}^\irr$ from the set of irreducible components of $\CV_\Delta$ to the set of irreducible components of $V_{\Delta,1}^\irr$, unless $(\Delta,M_0)\sim (\Delta_2,\ZZ^2)$, in which case, the surjection is from the subset $\{\CV_N\}_{N\ne N_0}$. 
    
    It remains to show that a general curve $[C]\in V_{\Delta,1}^\irr$ belongs to the image of a single irreducible component of $\CV_\Delta$. Let $(E;\bp)\in\CV_N$ and $(E';\bp')$ be two preimages of $[C]$. Then $E'=C^\nu=E$ and there exists a permutation $\sigma\in \prod_{j=1}^k{\mathfrak S}_{l_j}\subset {\mathfrak S}_l$ such that $p'_i=p_{\sigma(i)}$ for all $i$. However, the equations defining the loci $\CZ_N\subseteq \CZ_\Delta\subseteq \CM_{1,l}$ are invariant under the natural action of the subgroup $\prod_{j=1}^k{\mathfrak S}_{l_j}\subset {\mathfrak S}_l$ on the moduli space $\CM_{1,l}$, cf. Remark~\ref{rem:expCZ}. Thus, $\CZ_N$ is invariant under $\prod_{j=1}^k{\mathfrak S}_{l_j}$, and therefore so is $\CV_N$. Hence $(E';\bp')\in\CV_N$, and we are done.
\end{proof}

\subsection{The proof of the Main Theorem}

First, notice that condition (1) of the Main Theorem is equivalent to $M$ being an intermediate affine lattice $M_0\subseteq M\subseteq \ZZ^2$. 

Second, note that since there is a natural bijection between the intermediate affine sublattices $M_0\subseteq M\subseteq \ZZ^2$ and intermediate sublattice $N_0\subseteq N\subseteq \ZZ^2$, Corollaries \ref{cor:compCZ} and \ref{cor:compcorr} imply that the number of irreducible components of the Severi variety $V_{\Delta,1}^\irr$ is equal to the number of affine sublattices satisfying condition (1) of the Main Theorem unless one of the following holds: either (a) $\lw_{M_0}(\Delta)=1$ or (b) $(\Delta,M_0)\sim (\Delta_2,\ZZ^2)$, in which case, the number of irreducible components is one less than the number of such sublattices. Therefore, it remains to prove that an intermediate sublattice $M_0\subseteq M\subseteq \ZZ^2$ violates condition (2) of the Main Theorem if and only if $M=M_0$ and either (a) or (b) holds. 

Third, note that for any $M_0\subsetneq M\subseteq\ZZ^2$, we have $|\Delta^\circ\cap M|>|\Delta^\circ\cap M_0|\ge 0$ by Pick's formula, since $\partial\Delta\cap M=\partial\Delta\cap M_0$. Therefore, if $M$ violates condition (2), then $M=M_0$. It remains to prove the following elementary lemma:
\begin{lemma}
    Let $\Delta\subset\RR^2$ be a lattice polygon. Then $\Delta^\circ\cap\ZZ^2=\emptyset$ if and only if $\Delta$ is either a width-one polygon or twice a primitive triangle.
\end{lemma}

\begin{proof}
    The ``if'' part is obvious. In the opposite direction, assume that $\Delta^\circ\cap\ZZ^2=\emptyset$ and $\Delta$ is not a width-one polygon. Let us show that $\Delta$ is twice a primitive triangle. 
    
    Without loss of generality, $\partial_1\Delta$ is the longest facet of $\Delta$. By the width assumption, there exists $m\in \Delta\cap M$ such that the Euclidean area of $T:=\ConH\{\partial_1\Delta,m\}$ is at least $l_1$. Denote the sides of the triangle $T$ by $T_1:=\partial_1\Delta$, $T_2, T_3$, and their integral lengths by $l_{T_i}$. After changing the choice of $m$ if necessary, we may assume that $\max\{l_{T_2},l_{T_3}\}\le 2$; and without loss of generality $l_{T_2}\ge l_{T_3}$. Note that $T$ contains no inner integral points since $T\subseteq \Delta$. Thus, by Pick's formula, 
    \begin{equation}\label{eq:pickcalc}
        2l_1\le 2{\rm Area(T)}=l_1+l_{T_2}+l_{T_3}-2\le l_1+2l_{T_2}-2.
    \end{equation}
    Hence $l_{T_2}>1$, and therefore $l_{T_2}=2$. Since $\Delta$ is convex and $\Delta^\circ\cap\ZZ^2=\emptyset$, it follows that $T_2\subset\partial\Delta$. Thus, $l_1\ge l_{T_2}=2$ since $\partial_1\Delta$ is the longest facet of $\Delta$. Now, all inequalities in \eqref{eq:pickcalc} are necessarily equalities. Hence $l_1=l_{T_2}=l_{T_3}=2={\rm Area(T)}$, i.e., $T$ is twice a primitive triangle. Arguing as before, we see that $T_3\subset\partial\Delta$, which implies $\Delta=T$, and we are done.
\end{proof}

\bibliographystyle{amsalpha}
\bibliography{ref}

\providecommand{\bysame}{\leavevmode\hbox to3em{\hrulefill}\thinspace}
\providecommand{\MR}{\relax\ifhmode\unskip\space\fi MR }
\providecommand{\MRhref}[2]{%
  \href{http://www.ams.org/mathscinet-getitem?mr=#1}{#2}
}
\providecommand{\href}[2]{#2}
\begin{thebibliography}{CDGK23}

\bibitem[BLC23]{BL23}
Andrea Bruno and Margherita Lelli-Chiesa, \emph{Irreducibility of severi varieties on k3 surfaces}, 2023, arXiv:2112.09398.

\bibitem[Bou16]{Bou16}
David Bourqui, \emph{Algebraic points, non-anticanonical heights and the {S}everi problem on toric varieties}, Proc. Lond. Math. Soc. (3) \textbf{113} (2016), no.~4, 474--514. \MR{3556489}

\bibitem[CC99]{CC99}
Luca Chiantini and Ciro Ciliberto, \emph{On the {S}everi varieties on surfaces in {${\mathbb P}^3$}}, J. Algebraic Geom. \textbf{8} (1999), no.~1, 67--83.

\bibitem[CDGK23]{CDGK23}
Ciro Ciliberto, Thomas Dedieu, Concettina Galati, and Andreas~Leopold Knutsen, \emph{Nonemptiness of {S}everi varieties on {E}nriques surfaces}, Forum Math. Sigma \textbf{11} (2023), Paper No. e52, 32. \MR{4603111}

\bibitem[CFGK17]{CFGK17}
Ciro Ciliberto, Flaminio Flamini, Concettina Galati, and Andreas~Leopold Knutsen, \emph{Moduli of nodal curves on {$K3$} surfaces}, Adv. Math. \textbf{309} (2017), 624--654. \MR{3607287}

\bibitem[CGY23]{CGY23}
Nathan Chen, Fran\c{c}ois Greer, and Ruijie Yang, \emph{Nodal elliptic curves on {K}3 surfaces}, Math. Ann. \textbf{386} (2023), no.~3-4, 2349--2370. \MR{4612420}

\bibitem[CHT22]{CHT22}
Karl Christ, Xiang He, and Ilya Tyomkin, \emph{Degeneration of curves on some polarized toric surfaces}, J. Reine Angew. Math. \textbf{787} (2022), 197--240. \MR{4431908}

\bibitem[CHT23]{CHT23}
\bysame, \emph{On the {S}everi problem in arbitrary characteristic}, Publ. Math. Inst. Hautes \'{E}tudes Sci. \textbf{137} (2023), 1--45. \MR{4588594}

\bibitem[CHT25]{CHT25}
Karl Christ, Xiang He, and Ilya Tyomkin, \emph{The irreducibility of {H}urwitz spaces and {S}everi varieties on toric surfaces}, 2025, arXiv:2501.16238.

\bibitem[CLS11]{CLS11}
David~A. Cox, John~B. Little, and Henry~K. Schenck, \emph{Toric varieties}, Graduate Studies in Mathematics, vol. 124, American Mathematical Society, Providence, RI, 2011. \MR{2810322}

\bibitem[GLS07]{GLS07}
Gert-Martin Greuel, Christoph Lossen, and Eugenii~I. Shustin, \emph{Introduction to singularities and deformations}, Springer Monographs in Mathematics, Springer, Berlin, 2007. \MR{2290112}

\bibitem[Har86]{Har86}
Joe Harris, \emph{On the {S}everi problem}, Invent. Math. \textbf{84} (1986), no.~3, 445--461. \MR{837522 (87f:14012)}

\bibitem[KS13]{KS13}
Steven~L. Kleiman and Vivek~V. Shende, \emph{On the {G}\"{o}ttsche threshold}, A celebration of algebraic geometry, Clay Math. Proc., vol.~18, Amer. Math. Soc., Providence, RI, 2013, With an appendix by Ilya Tyomkin, pp.~429--449.

\bibitem[Lan20]{L20}
Lionel Lang, \emph{Monodromy of rational curves on toric surfaces}, J. Topol. \textbf{13} (2020), no.~4, 1658--1681. \MR{4186140}

\bibitem[LT23]{LT23}
Lionel Lang and Ilya Tyomkin, \emph{A note on the {S}everi problem for toric surfaces}, Math. Ann. \textbf{385} (2023), no.~3-4, 1677--1705. \MR{4566703}

\bibitem[Pic99]{Pic99}
Georg~A. Pick, \emph{Geometrisches zur {Zahlenlehre}}, Sonderabdr. {Naturw}.-medizin. {Verein} f. {B{\"o}hmen} ``{Lotos}'' {Nr}. 8, 9 {S}. {{\(8^{\circ}\)}} (1899), 1899.

\bibitem[Smi61]{Smi61}
Henry J.~Stephen Smith, \emph{On systems of linear indeterminate equations and congruences}, Philosophical Transactions of the Royal Society of London \textbf{151} (1861), 293--326.

\bibitem[Tes09]{Tes09}
Damiano Testa, \emph{The irreducibility of the spaces of rational curves on del {P}ezzo surfaces}, J. Algebraic Geom. \textbf{18} (2009), no.~1, 37--61. \MR{2448278}

\bibitem[Tyo07]{Tyo07}
Ilya Tyomkin, \emph{On {S}everi varieties on {H}irzebruch surfaces}, Int. Math. Res. Not. IMRN (2007), no.~23, Art. ID rnm109, 31. \MR{2380003}

\bibitem[Tyo14]{Tyo14}
\bysame, \emph{An example of a reducible {S}everi variety}, Proceedings of the {G}\"{o}kova {G}eometry-{T}opology {C}onference 2013, G\"{o}kova Geometry/Topology Conference (GGT), G\"{o}kova, 2014, pp.~33--40.

\bibitem[Tyo24]{Tyo24}
Ilya Tyomkin, \emph{The geometry of {S}everi varieties}, 2024, arXiv:2411.11431.

\bibitem[Zah22]{Z22}
Adrian Zahariuc, \emph{The {S}everi problem for abelian surfaces in the primitive case}, J. Math. Pures Appl. (9) \textbf{158} (2022), 320--349. \MR{4360371}

\end{thebibliography}
\end{document}